\documentstyle[12pt]{article}

\newcommand{\commentout}[1]{}
\newcommand{\bes}{\begin{displaymath}}
\newcommand{\ees}{\end{displaymath}}
\newcommand{\be}{\begin{equation}}
\newcommand{\ee}{\end{equation}}
\newcommand{\bi}{{\bf i}}
\newcommand{\ba}{\begin{eqnarray}}
\newcommand{\ea}{\end{eqnarray}}
\newcommand{\bas}{\begin{eqnarray*}}
\newcommand{\eas}{\end{eqnarray*}}

\newcommand{\al}{\alpha}
\newcommand{\bE}{{\bf E}}
\newcommand{\bone}{{\bf 1}}

\newcommand{\bI}{{\bf I}}
\newcommand{\bW}{{\bf W}}
\newcommand{\bD}{{\bf D}}

\newcommand{\cL}{{\cal L}}
\newcommand{\cF}{{\cal F}}

\newcommand{\cG}{{\cal G}}

\newcommand{\bt}{\beta}

\newcommand{\bbm}{{\bf m}}
\newcommand{\bbsig}{\mbox{\boldmath{$\sigma$}}}
\newcommand{\bbi}{{\bf i}}
\newcommand{\bbl}{{\bf \ell}}
\newcommand{\tk}{\tilde{k}}
\newcommand{\si}{\sigma}

\newcommand{\vV}{{\bf V}}

\newcommand{\vphi}{\varphi}

\newcommand{\ep}{\varepsilon}
\newcommand{\vep}{\varepsilon}

\newcommand{\la}{\lambda}

\newcommand{\Om}{\Omega}

\newcommand{\ga}{\gamma}
\newcommand{\bR}{{\bf R}}

\newcommand{\bk}{{\bf k}}

\newcommand{\bze}{{\bf 0}}

\newcommand{\kropa}{$_{\Box}$}
\newcommand{\bU}{{\bf U}}

\newcommand{\da}{\downarrow0}
\newtheorem{theorem}{Theorem}
\newtheorem{proposition}{Proposition}

\newtheorem{lemma}{Lemma}
\newcommand{\nwc}{\newcommand}
\nwc{\m}{\mbox}
\nwc{\bm}{\boldmath}
\nwc{\ubm}{\unboldmath}
\nwc{\bmu}{\m{\bm $u$\ubm}}
\nwc{\bmx}{\m{\bm $x$\ubm}}
\nwc{\bx}{\bmx}

\nwc{\bmv}{\m{\bm $v$\ubm}}
\nwc{\cE}{{\cal E}}

\newcommand{\bbxo}{\bx_\ep(s)}

\nwc{\beq}{\begin{eqnarray}}
\nwc{\eeq}{\end{eqnarray}}
\begin{document}
\title{
Taylor-Kubo Formula
for Turbulent Diffusion in a Non-Mixing Flow with Long-Range Correlation}
\author{Albert C. Fannjiang \\ Department of
Mathematics, University of California, Davis\and Tomasz Komorowski \\
Department of Statistics, University of California, Berkeley}
\thanks{The research of A. F. is supported by National Science Foundation
Grant No. DMS-9600119.}

\maketitle

{\bf Abstract} We prove the Taylor-Kubo formula
 for a class of isotropic, non-mixing flows with long-range correlation.
 For the proof, we develop
 the method of high order correctors expansion.

 {\bf Keywords} Turbulent diffusion, mixing, corrector.

 {\bf AMS subject classification} Primary 60F05, 76F05, 76R50; Secondary 58F25.

 {\bf Abbreviated title} Taylor-Kubo formula for non-mixing flows.

\newpage

\section{Introduction}
Let the deterministic motion of a passive particle in a random velocity
field $\vV(t,\bmx)$ be described by
\be
\label{a1.1}
{d\bmx_\ep(t)\over dt}={1\over\ep}\vV({t\over \ep^2},\bmx_\ep(t))
\ee
where $\ep$ is a small parameter and
$\vV(t,\bmx)=(V_1(t,\bmx),\cdots,V_d(t,\bmx))$ a time-stationary,
space-homogeneous,
incompressible, $\nabla\cdot\vV=0$, zero-mean velocity field,
The question is to determine the asymptotics, as $\ep$ tends to zero,
of the displacement
$\bmx_\ep(t)$ from the statistics of the velocity.

Scaling like (\ref{a1.1}) arises naturally when the velocity
has a non-zero mean drift and a small fluctuation
$\bar{\vV}+\ep\vV'(t,\bmx)$. By changing to the moving coordinate
system $\bmx\to\bmx+\bar{\vV}t$ and taking the long time
limit $t\to t/\ep^2$ we have a equation of the form (\ref{a1.1})
with $\vV(t,\bmx)=\vV'(t,\bmx+\bar{\vV}t)$.

Under the assumption of {\em strong mixing in time} among other things,
Kesten and Papanicolaou \cite{KP} proved that $\bmx_\ep(t)$ converges
to a Brownian motion with the diffusion coefficient given by
a Taylor-Kubo formula
\be
\label{TK}
D_{ij}^*=\int\limits^\infty_0\left\{\bE[V_i(t,\bze)V_j(0,\bze)]+ \bE[V_i(t,\bze)V_j(0,\bze)]\right\}dt
\ee
(see also Carmona and Fouque \cite{CF} for the corresponding
result of stochastic flows). No assumption on the decay of
velocity decorrelation in space is required.

Convergence to Brownian motion in the absence of molecular diffusion
is often referred to as {\em turbulent diffusion}
because the large scale diffusive motion is a result of  random flows
instead of molecular diffusion.

\commentout{
Let us be more specific about slowly varying velocity
and the associated turbulent diffusion theorem. The
velocity field $\vV$ in (\ref{a1.0}) is assumed to have zero expectation
\[
\bE\{\vV\}=0.
\]
Let the velocity field $\vV$ in (\ref{a1.0}) be replaced
with a new velocity field $\vV_\ep$
\be
\label{a1.4}
\vV_\ep(t,\bmx)=\vV(t,\ep\bmx),
\ee
with a small parameter $\ep$,
to model slowly varying fields.
Let $\bmx_\ep$ be  the rescaled displacement
\be
\label{a1.2}
\bmx_\ep(t)=\ep\bmx({t\over \ep^2}).
\ee
The scaling in (\ref{a1.2}) is called  {\em diffusive}.
}

Does the turbulent diffusion
theorem hold for velocities
lacking the temporal mixing property but with correlations decaying in space?
To study the
interplay of temporal and spatial correlations in a precise way,  we consider
the class
of time-stationary, space-homogeneous, isotropic, Markov, Gaussian velocity
fields $\vV(t,\bmx)$, defined
on the probability space
$(\Om,{\cal V}, P)$, with two-point correlation tensor $\bR=[R_{ij}(t,\bmx)]$
 given by the  Fourier transform
\be
\label{a1.3}
R_{ij}(t,\bmx)=\bE\left[\vV(t,\bmx)\otimes\vV(0,\bze)\right]=
\ee
\[
\int_{R^d}\cos{(\bk\cdot\bx)}e^{-|\bk|^{2\bt}t}\hat{\bR}(\bk)d\bk
\]
where
$\bE$ stands for the expectation.
The spectral density $\hat{\bR}(\bk)$ is given by
a power law
\be
\label{060103}
\label{a1.3b}
\hat{\bR}(\bk)=\frac{a(|\bk|)}{|\bk|^{2\al+d-2}}\left(\bI-\frac{\bk\otimes\bk}{|\bk|^2}\right),\quad\alpha<1,\quad\beta\geq 0,
\ee
with a compactly supported,
continuous function
$a:[0,+\infty)\rightarrow R_+$.
The factor $\bI-\frac{\bk\otimes\bk}{|\bk|^2}$ in (\ref{a1.3b}) ensures
that the velocity
field  is incompressible.

The function $\exp{(-|\bk|^{2\beta}t)}$ in (\ref{a1.3}) is called the
{\em time correlation function} of the velocity $\vV$.
The  spectral density $\hat{\bR}(\bk)$ is integrable over $\bk$
 for $\alpha<1$ and, thus, (\ref{a1.3})-(\ref{060103}) defines
a velocity field with a finite second moment. The ultraviolet
cut-off $K$ is needed to avoid divergence of integral over large $|\bk|$.
Because of the ultraviolet cutoff $\vV$
is jointly continuous in $t,\bmx$ and is $C^\infty$ in $\bmx$ almost surely.

The parameter $\alpha$ is directly
related to the decay exponent of $\bR: \bR(\bmx)\sim |\bmx|^{\al-1}$
for $|\bmx|\gg1$.
As $\alpha$ increases to one, the spatial decay exponent of $\bR$
decreases to zero and, consequently, spatial correlation of velocity
increses.
On the other hand,
for $\beta>0$, the velocity field
lacks the spectral gap and thus strong mixing property (see \cite{R}).
We restrict our attention to the case $\alpha<1,\beta >0$ which
corresponds to a velocity field
with arbitrarily long scales but not
the strong mixing property.

What is the region in the $(\alpha, \beta)$ plane where the
turbulent diffusion
theorem,
with the Taylor-Kubo formula, holds? It is easy to find the necessary
 condition by imposing the convergence of
the Taylor-Kubo
formula. A plain calculation
\beq
\label{0220}
D^*_{ij}&=&\int\limits^\infty_0 R_{ij}(t,0)dt\\
&=&\int\limits_{R^d}(\delta_{ij}-{k_ik_j\over |\bk|^2})
a(|\bk|)|\bk|^{2-2\alpha-d}
\int\limits^\infty_0 \exp{(-|\bk|^{2\beta}t)}dt d\bk \label{a1.10}\\
&=&\int\limits_{R^d}(\delta_{ij}-{k_ik_j\over |\bk|^2})
{a(|\bk|)\over |\bk|^{2\alpha+2\beta-1}}{d\bk\over |\bk|^{d-1}}
\eeq
leads to
the condition
\be
\label{a1.5}
\alpha+\beta<1.
\ee
As it turns out, (\ref{a1.5}) is also sufficient.
This is our main result
which generalizes the turbulent diffusion theorem of \cite{KP}
to a class of non-mixing flows.
\begin{theorem}
\label{theorem1}
Let $\vV(t,\bx)$ be
a Markov, Gaussian velocity
field with correlation given by (\ref{a1.3})-(\ref{a1.3b}).
Then,
for $\al+\bt<1$, the invariance principle holds
for the displacement $\bmx_\ep(t)$. Namely, $\bmx_\ep(t)$,
as continuous
process,
converges weakly, as $\ep\da$,
to a Brownian Motion with
covariance matrix given by formula (\ref{TK}).
\end{theorem}
Additional molecular diffusion in (\ref{a1.1})
would act like a regular perturbation
to the eddy diffusivity defined by the Taylor-Kubo formula.


What happens if (\ref{a1.5}) is violated?
From
the divergence of the Taylor-Kubo formula one judges
that the long time asymptotics
should be faster than diffusion and, consequently, a
different scaling limit is necessary.
It turns out that under an anomalous scaling the
displacement converges
to a fractional Brownian motion. The fractional Brownian motion
limit theorem will be studied in a forthcoming paper.
The dichotomy is conveniently represented in the figure.

\section{Multiple stochastic integrals}

By the Spectral Theorem (see,  e.g., \cite{rozanov}) we assume
 without loss of generality that there exist two
independent, identically distributed, real vector valued, Gaussian spectral measures
$\hat{\vV}_{l}(t,\cdot)$, $l=0,1$ such that
\be
\vV(t,\bx)=\int \hat{\vV}_0(t,\bx,d\bk),
\label{a.0}
\ee
where
\[
\hat{\vV}_0(t,\bx,d\bk):=c_0(\bk\cdot\bx)\hat{\vV}_0(t,d\bk)+c_1(\bk\cdot\bx)\hat{\vV}_1(t,d\bk)
\]
with $c_0(\phi)\equiv\cos{(\phi)},c_1(\phi)\equiv\sin{(\phi)}$.
Define also
\[
\hat{\vV}_1(t,\bx,d\bk):=-c_1(\bk\cdot\bx)\hat{\vV}_0(t,d\bk)+c_0(\bk\cdot\bx)\hat{\vV}_1(t,d\bk).
\]
We have the relation
\beq
\label{a.2}
\partial \hat{\vV}_0(t,\bx,d\bk)/\partial x_j&= &k_j \hat{\vV}_1(t,\bx,d\bk),\\
\quad\partial \hat{\vV}_1(t,\bx,d\bk)/\partial x_j&=& -k_j \hat{\vV}_0(t,\bx,d\bk).\label{a.1}
\eeq
One can check  that
$\int \hat{\vV}_1(t,\bx,d\bk)$
is a random field distributed identically to and independently of $\vV$.


We define the {\em multiple stochastic integral}
\be
\label{08144}
\int\cdots\int  \psi(\bk_1,\cdots,\bk_N)\widehat{\vV}_{l_1}(t_1,\bx_1,d\bk_1)\otimes\cdots\otimes\widehat{\vV}_{l_N}(t_N,\bx_N,d\bk_N)
\ee
for any  $l_1,\cdots,l_N\in\{0,1\}$ and a suitable family of functions $\psi$
by using  the {\em Fubini theorem} (see (\ref{061902}) below).
For $\psi_1,\cdots,\psi_N\in {\cal S}(R^d)$,
the Schwartz space, and $l_1,\cdots,l_N\in\{0,1\}$  we set
\beq
\label{061902}
&&
\int\cdots\int \psi_1
(\bk_1)\cdots\psi_N(\bk_N)\widehat{\vV}_{l_1}(t_1,\bx_1,d\bk_1)\otimes\cdots\otimes\widehat{\vV}_{l_N}(t_N,\bx_N,d\bk_N)\\
&:=&
\int\psi_{1}(\bk_1)\widehat{\vV}_{l_1}(t_1,\bx_1,d\bk_1)\otimes\cdots\otimes\int\psi_N(\bk_N)\widehat{\vV}_{l_N}(t_N,\bx_N,d\bk_N).\nonumber
\eeq
We then extend the definition of multiple integration to the closure ${\cal H}$ of the Schwartz space ${\cal S}((R^d)^N,R)$ under the norm
\be
\label{0807}
\|\psi\|^2:=\int\cdots\int \psi(\bk_1,\cdots,\bk_N) \psi(\bk_1',\cdots,\bk_N')
\ee
\[
\bE\left[\widehat{\vV}_{l_1}(t_1,\bx_1,d\bk_1)\otimes\cdots\otimes\widehat{\vV}_{l_N}(t_N,\bx_N,d\bk_N):\widehat{\vV}_{l_1}(t_1,\bx_1,d\bk_1')\otimes\cdots\otimes\widehat{\vV}_{l_N}(t_N,\bx_N,d\bk_N')\right].
\]
The expectation is to be calculated by the formal rule
\[
\bE\left[\widehat{V}_{l,i}(t,\bx,d\bk)\widehat{V}_{l',i'}(t',\bx',d\bk')\right]=e^{-|\bk|^{2\bt}|t-t'|}\delta_{l,l'}c_0(\bk\cdot(\bx-\bx'))
\widehat{R}_{i,i'}(\bk)\delta(\bk-\bk')d\bk d\bk'.
\]
Among various possible definitions of multiple stochastic integral (see, e.g.,
\cite{rozanov2,shiryaev,sinai}) this approach (\cite{shiryaev})
seems to best serve our purposes.


\label{061901}

We denote the stochastic integral
(\ref{08144}) by $\Psi_{\bbl}(t_1,\cdots,t_N,\bx_1,\cdots,\bx_N)$ or,
sometimes, simply by $\Psi_{\bbl}$. Here
$\bbl=(l_1,\cdots,l_N), l_1,\cdots,l_N\in
\{0,1\}$. Sometimes when  $\bbi=(i_1,\cdots,i_d)$, with $i_1,\cdots,i_d\in
\{1,2,\cdots,d\}$ is fixed we shall denote the corresponding component of
the above tensor valued random field by $\Psi_{\bbl,\bbi}$.

Note that $\Psi_{\bbl,\bbi}\in
H^N(\vV)$ - the Hilbert space
obtained as a completion of the space of $N$-th degree polynomials in
variables $\int\psi(\bk)\widehat{\vV}(t,\bx,\bk)$ with
respect to the standard $L^2$ norm.
\begin{proposition}
\label{proposition611}
For any $(t_1,\bx_1),\cdots,(t_N,\bx_N)\in R\times R^d$ and $p>0$,
$\Psi_{\bbl,\bbi}$ belongs to  $L^p(\Om)$ and
\be
\label{061701}
\left(\bE|\Psi_{\bbl,\bbi} |^p\right)^{1/p}\leq C\left(\bE|\Psi_{\bbl,\bbi}|^2\right)^{1/2}
\ee
with the constant $C$ depending only on $p,N$ and the dimension $d$.
Moreover, $\Psi_{\bbl,\bbi}$ is differentiable in the mean square sense with
\be
\label{061702}
\nabla_{\bx_j}\Psi_{\bbl,\bbi}(t_1,\cdots,t_N,\bx_1,\cdots,\bx_N)=
\int\cdots\int \bk_j \psi(\bk_1,\cdots,\bk_N)
\ee
\[
\widehat{V}_{l_1,i_1}(t_1,\bx_1,d\bk_1)\cdots\widehat{V}_{1-l_j,i_j}(t_j,\bx_j,d\bk_j)\cdots\widehat{V}_{l_N,i_N}(t_N,\bx_N,d\bk_N).
\]
\end{proposition}

The proof of Proposition
\ref{proposition611} is standard and follows directly from the well known
hypercontractivity property for Gaussian measures
(see, e.g., \cite{J}, Theorem 5.1. and its corollaries), so
we do not repeat it here.

The field $\vV$ is Markovian
\be
\label{1.10}
\bE\left[\int\psi(\bk)\widehat{\vV}_l(t,\bx,d\bk)\left|\right.{\cal
    V}_{-\infty,s}\right]= \int
e^{-|\bk|^{2\bt}(t-s)}\psi(\bk)\widehat{\vV}_l(s,\bx,d\bk),~~l=0,1,
\ee
 for  all
 $\psi\in{\cal S}(R^d,R)$,
where ${\cal
    V}_{a,b}$ denotes the $\si$-algebra generated by random variables
    $\vV(t,\bx)$, for $t\in [a,b]$ and $\bx\in R^d$.

To calculate expectation of multiple product of Gaussian random variables, it is
convenient to use a graphical representation, borrowed from quantum
field theory. We refer to, e.g., Glimm and Jaffe \cite{GJ} and
Janson \cite{J}.
A {\em
Feynman diagram} ${\cal F}$ (of order $n\geq 0$ and rank $r\geq0$)
is a graph consisting of a set $B(\cF)$ of $n$ vertices and a set $E(\cF)$
of $r$
edges without common endpoints. So there are $r$ pairs
of vertices, each joined by an edge, and $n-2r$ unpaired vertices,
called  free vertices. $B(\cF)$ is a set
of positive integers.
An edge whose endpoints are $m,n\in B$ is represented by $\widehat{mn}$
(unless otherwise specified, we always assume $m<n$);
and an edge includes its endpoints.
A diagram ${\cal F}$ is said
to be {\em based on} $B(\cF)$.
Denote the set of
free vertices (i.e. points which are not
endpoints of any edges) by $A({\cal F})$, so $A({\cal F})=\cF\setminus
E(\cF)$.
The
diagram is {\em complete} if $A({\cal F})$ is empty and
{\em incomplete}, otherwise.
Denote by $\cG(B)$ the set of all diagrams based on $B$,
by $\cG_c(B)$ the set of all complete diagrams based on $B$ and
by $\cG_i(B)$ the set of all incomplete diagrams based on $B$.
A special class of  diagrams, denoted by
$\cG_s(B)$, plays an important role  in the subsequent
analysis: a diagram $\cF$ of order $n$ belongs to $\cG_s(B)$
if $A_k(\cF)$ is not empty for all $k=0,1,...,n.$
Let $B=\{1,2,3,...,n\}$. Denote by ${\cal F}_{|k}$ the sub-diagram of ${\cF}$,
based on $\{1,\cdots,k\}$.
Define $A_k({\cal F})=A({\cal F}_{|k})$.

We work out the conditional expectation for
multiple spectral integrals using the Markov property (\ref{1.10}).
\begin{proposition}
For any function $\psi\in {\cal H}$ and $l_1,\cdots,l_N\in \{0,1\}$,
$i_1,\cdots,i_N\in\{1,\cdots,d\}$,
\label{prop1}
\be
\label{060106}
\bE\left[\int\cdots\int
\psi(\bk_1,\cdots,\bk_N)
\widehat{V}_{l_1,i_1}(t,\bx_1,d\bk_1)\cdots
\widehat{V}_{l_N,i_N}(t,\bx_N,d\bk_N)\left|\right.{\cal V}_{-\infty,s}\right]=
\ee
\[
\sum\limits_{{\cal
   F}\in\cG(\{1,...,N\})}\int\cdots\int\exp\left\{-\sum\limits_{m\in A({\cal
   F})}|\bk_m|^{2\bt}(t-s)\right\}\psi(\bk_1,\cdots,\bk_N)
\widehat{\vV}_{s,\bx_1,\cdots,\bx_N}(d\bk_1,\cdots,d\bk_N;{\cal
    F})
\]
with
\be
\label{01201}
\widehat{V}_{s,\bx_1,\cdots,\bx_N}(d\bk_1,\cdots,d\bk_N;{\cal
   F}):=\prod\limits_{m\in A({\cal
  F})}\widehat{V}_{l_m,i_m}(s,\bx_m,d\bk_m)
\ee
\[
\mathop{\prod\limits_{\widehat{mn}\in
   E({\cal F})}}
\left[1-e^{-\left(|\bk_{m}|^{2\bt}+|\bk_{n}|^{2\bt}\right)(t-s)}
\right]\bE[\widehat{V}_{l_{m},i_m}(s,\bx_m,d\bk_{m})
\widehat{V}_{l_{n},i_n}(s,\bx_n,d\bk_{n})].
\]
\end{proposition}

{\bf Proof.} Without loss of generality we consider
$
\psi(\bk_1,\cdots,\bk_N)=\bone_{A_1}(\bk_1)\cdots\bone_{A_N}(\bk_N)
$
for some Borel sets $A_1,\cdots,A_N$.

$\widehat{\vV}_l(t,A_i)=\widehat{\vV}^{0}_l(t,A_i)+\widehat{\vV}^{1}_l(t,A_i)$ where
$\widehat{\vV}^{0}_l(t,\cdot)$, $\widehat{\vV}^{1}_l(t,\cdot)$ are
the orthogonal projection of $\widehat{\vV}_l(t,\cdot)$ on the space
$L^2_{-\infty,s}$ and its complement, respectively.
Here  $L^2_{a,b}$ denotes $L^2$
closure of the linear span over $\vV(s,\bx)$, $a\leq s\leq b$, $\bx\in R^d$.
The conditional expectation in
(\ref{060106}) equals
\[
\sum\limits_{{\cal
    F}\in\cG(\{1,...,N\})}\mathop{\prod\limits_{\widehat{mn}\in
       E({\cal F})}}
     \bE\left[\widehat{V}_{i_{m},l_{m}}^1(t,A_{m})\widehat{V}_{i_{n},l_{n}}^1(t,A_{n})\right]\prod\limits_{m\in A({\cal F})}\widehat{V}_{i_m,l_m}^0(t,A_{m}).
\]
The statement follows upon the application of the relations
\[
\widehat{\vV}_{l}^0(t,A)=\int\limits_{A}e^{-|\bk|^{2\bt}(t-s)}\widehat{\vV}_{l}(s,d\bk)
\]
and
\[
\bE\left[\widehat{\vV}_{l}^1(t,A)\otimes\widehat{\vV}_{l'}^1(t,B)\right]=
\]
\[
\int\limits_{A}\int\limits_{B}\delta_{l,l'}\left\{\bE\left[\widehat{\vV}_{l}(t,d\bk)\otimes\widehat{\vV}_{l'}(t,d\bk')\right]-\bE\left[\widehat{\vV}_{l}^0(t,d\bk)\otimes\widehat{\vV}_{l'}^0(t,d\bk')\right]\right\}.
\]
\quad\rule{2mm}{3mm}

\label{sec2}
For fixed $\bx\in R^d$ we define a process
\[
\Psi_{\bbl}(t;\bx):=
\Psi_{\bbl}(t,\cdots,t,\bx,\cdots,\bx),\qquad\mbox{for }t\geq0.
\]
Let us introduce, after Chapter 3 of \cite{kushner},
the {\em pseudogenerator} of the
above process ${\cal L}(t;\bx)$, $t\geq0$ as
\[
\lim\limits_{\delta\da}\frac{\bE\left[\Psi_{\bbl}(t+\delta;\bx)\left|\right.{\cal
  V}_{-\infty,t}\right]-\Psi_{\bbl}(t;\bx)}{\delta}
\]
in the $L^1$ sense. After elementary calculations using (\ref{060106})
we obtain that
\be
\label{01202}
{\cal L}(t;\bx)=\sum\limits_{\widehat{pq}}
\int\cdots\int(|\bk_p|^{2\bt}+|\bk_q|^{2\bt})
\psi(\bk_1,\cdots,\bk_N)\delta(\bk_p-\bk_q)\delta_{l_p,l_q}
\widehat{\vV}_{\widehat{pq}}(s,\bx,d\bk_1,\cdots,d\bk_N)
d\bk_pd\bk_q.
\ee
Here the summation extends over all edges $\widehat{pq}$ with
vertices belonging to
the set ${\cal S}=\{1,\cdots,N\}$ and
$\widehat{\vV}_{\widehat{pq}}
(s,\bx,d\bk_1,\cdots,d\bk_N)$ indicates that the multiple
stochastic integration
is taken over the product of stochastic measures
$\widehat{\vV}_{l_n}(s,\bx,d\bk_n)$ for all $n\in{\cal S}\setminus\{p,q\}$.

The fundamental  property
 of ${\cal L}(t;\bx)$, $t\geq0$ is that for any fixed $\bx$
\be
\label{01181}
M(t;\bx):=\Psi_{\bbl}(t;\bx)-\int\limits_0^t{\cal L}(s;\bx)ds
\ee
is a continuous trajectory martingale, cf. \cite{kushner}. In fact  we can
give an explicit expression on its quadratic variation process $<M>_t$,
$t\geq0$ in terms of
the  spectral representation.
\begin{proposition}
\label{prop5}
For any fixed $\bx_1,\cdots,\bx_N\in R^d$
\be
\label{01211}
<M(,\cdot;\bx)>_t=\int\limits_0^tW(s,\bx)ds
\ee
where
\[
W(t,\bx):=2\sum\limits_{\widehat{pq}}\int\cdots\int\int\cdots\int\psi(\bk_1,\cdots,\bk_N)
\psi(\bk_1',\cdots,\bk_N')a(|\bk_p|)|\bk_p|^{2\bt+2-2\al-d}
\]
\[
\left(\bI-
\frac{\bk_p\otimes\bk_p}{|\bk_p|^2}\right)\otimes\widehat{\vV}_{\widehat{p,q}}(t,\bx,d\bk_1,\cdots,d\bk_{2N})\delta(\bk_p-\bk_q)
\delta_{l_p,l_q}d\bk_p d\bk_q
d\bk_pd\bk_q.
\]
Here we adopt the convention that $\bk_n'=\bk_{N+n}$, $l_n=l_{N+n}$.
The summation extends over all edges $\widehat{pq}$ with
vertices $p,q$ belonging to the sets $\{1,\cdots,N\}$ and
$\{N+1,\cdots,2N\}$ respectively.
\end{proposition}
{\bf Proof.} The quadratic variation of $M(\cdot;\bx)$ coincides with the pseudogenerator
of $M(\cdot;\bx)\otimes M(\cdot;\bx)$. After an elementary calculation we get that the left hand side
of (\ref{01211}) equals
\[
{\cal L}'(t;\bx)-
\Psi_{\bbl}(t;\bx)\otimes{\cal L}(t;\bx)-{\cal L}(t;\bx)\otimes
\Psi_{\bbl}(t;\bx),
\]
with ${\cal L}'(t;\bx)$, $t\geq0$ denoting the pseudogenerator of
$\Psi_{\bbl}(t;\bx)\otimes\Psi_{\bbl}(t;\bx)$, $t\geq0$.
The proof of the proposition can be then concluded from an
application of formula (\ref{01202})\quad\rule{2mm}{3mm}
\section{Proof of Theorem \ref{theorem1}}
We begin this section by defining the key concepts:
{\em the $n$-th
$\la$-corrector},
$\chi_{\la}^{(n)}$, and {\em the $n$-th $\la$-convector},
$\bU^{(n)}_{\la}$ recursively:
\beq
\bU^{(0)}_{\la}(s,\bx)&=&\vV(s,\bx)\\
\label{123}
\chi_{\la}^{(n)}(t,\bx)&:=&\int\limits_t^{+\infty}e^{-\lambda(s-t)}
\bE\left[\bU^{(n-1)}_{\la}(s,\bx)\left|\right.{\cal
  V}_{-\infty,t}\right]ds\\
\label{122}
\bU^{(n)}_{\la}(t,\bx)&:=&\vV(s,\bx)\cdot\nabla\chi_{\la}^{(n)}(t,\bx).
\eeq
The differentiation in (\ref{122}) is taken in the
mean square sense.

The definitions requires the integrability of
convectors and the differentiability
of correctors,
as stated in the following lemma.
\begin{lemma}
\label{lemma2}
 $\chi_{\la}^{(n)}$, $\bU^{(n)}_{\la}$
are  $L^p$ integrable for any $p\geq 1$,
$\la>0$ and $n\geq0$ and
satisfy the following hypercontractivity property
\be
\label{061201}
\left(\bE|\chi_{\la}^{(n)}(t,\bx)|^p\right)^{1/p}\leq C\left(\bE|\chi_{\la}^{(n)}(t,\bx)|^2\right)^{1/2}
\ee
\[
\left(\bE|\bU^{(n)}_{\la}(t,\bx)|^p\right)^{1/p}\leq C\left(\bE|\bU^{(n)}_{\la}(t,\bx)|^2\right)^{1/2},
\]
where the constant $C>0$ depends on $n,p,d$ but not on $\la$.

There exist versions of  $\chi_{\la}^{(n)}$ and $\bU_{\la}^{(n)}$
which are jointly stationary in the strict sense,  continuous in both $(t,\bx)$, $C^\infty$ smooth in
$\bx$.
\end{lemma}

$L^p$ integrability follows from  the
hypercontractivity property on  Gaussian Measure spaces
(see Theorem 5.1, Nelson's theorem and the corollaries in \cite{J}).
The existence of  regular versions of the fields  can be proved by
using aforementioned Nelson's theorem in conjunction with Kolmogorov's
criterion of continuous processes
(see Theorem 3.2.5 and corollary, \cite{adler}).

This and other lemmas in
this section are proved in the next section.

We define {\em the $n$-th rescaled corrector along the path} by
\be
\label{126}
\chi_{\ep}^{(n)}(t):=\ep^n\chi_{\ep^2}^{(n)}(\frac{t}{\ep^{2}},\bx_\ep(t))
\ee
and {\em the  $n$-th rescaled convector} by
\be
\label{129}
\bU_{\ep}^{(n)}(t):=\ep^{n-1}\bU_{\ep^2}^{(n)}(\frac{t}{\ep^{2}},\bx_\ep(t)).
\ee
Because of the incompressibility of $\vV$,  both
$\chi_{\ep}^{(n)}(t)$ and
$\bU^{(n)}_{\ep}(t)$ are stationary (see
Theorem 2 of  \cite{port-stone}).

Next we represent $\chi_{\la}^{(n)}$ by
multiple spectral integrals.
To simplify the presentation we introduce some notations.
For $\bbm=(m_1,\cdots,m_p)$ and
$p'\leq p$ we
set $\bbm_{|p'}:=(m_1,\cdots,m_{p'})$, $\bbm^{|p'}:=(m_{p'},\cdots,m_{p})$ and
$(\bbm, k):=(m_1,\cdots,m_p,k)$.

 \commentout{
Let ${\cal Z}_n$ be the family of
  positive-integer-valued sequences $\bbm=(m_1,\cdots,m_n)$  of length $n$
   such that $m_1:=1$ and $m_l<l\leq n$.
   By induction we define functions ${\cal S}_n$ from
   $\bbm\in {\cal Z}_n$ to $\bbsig=(\sigma_1,\cdots,\sigma_n)\in\{0,1\}^n$:
   ${\cal S}_1((1)):=(0)$. Suppose
   ${\cal S}_n(\bbm)=\bbsig$. Then
   \be
   {\cal S}_{n+1}((\bbm, p))=\bbsig'= (\si_1',\cdots,\si_{n+1}')
   \ee
   with
   \beq
   \si_1'&:=&0\\
   \si_{p+1}'&:=&1-\si_p.\\
   \si_{j+1}'&:=&\si_j\quad\mbox{otherwise}.
   \eeq

 \commentout{
Let ${\cal Z}_n$ be the family of
 positive-integer-valued sequences $\bbm=(m_1,\cdots,m_n)$  of length $n$
 such that $m_1:=1$ and $m_l<l, \forall l\leq n$.
By induction we define functions ${\cal S}_n(\bbm)=\bbsig$ from
$\bbm\in {\cal Z}_n$ to $\bbsig=(\sigma_1,\cdots,\sigma_n)\in\{0,1\}^n$
\beq
\label{070805}
{\cal S}_1((1))&:=&(0)\\
\label{070806}
{\cal S}_{n+1}((\bbm, p))&=&(\si_1',\cdots,\si_{n+1}'),\quad\hbox{for  }p<n+1
\eeq
with
\beq
\nonumber
\si_1'&:=&0\\
\si_{j+1}'&:=&\si_j\quad\mbox{for}\,\,1\leq j\leq n,\quad j\not=p
\nonumber\\
\si_{p+1}'&:=&1-\si_p.
\eeq
}}

We introduce next the {\em proper} functions of order $n,
 \bbsig:\{1,\cdots,n\}\rightarrow \{0,1\}$ that appears in the
 statement of the next lemma. The proper function of order 1
 is unique and is given
 by $\bbsig(1)=0$. Any proper function, $\bbsig'$, of order $n+1$ is generated
 from a proper function $\bbsig$ of order $n$ as follows.
 For some $p\leq n$,
 \beq
 \label{299}
 \bbsig'(n+1)&:=&0\\
 \nonumber
 \bbsig'(k)&:=&\bbsig(k)\quad
 \mbox{for}\quad k\leq n\mbox{ and }k\not=p\\
 \bbsig'(p)&:=&1-\bbsig(p).
 \eeq
 In other words, each proper function $\bbsig$ of order
 $n$ generates $n$ different proper functions of order
 $n+1$. Thus, the total number of proper functions of order $n$
 is $(n-2)\!$.

\begin{lemma}
\label{lemma2b}
For any $i_1\in\{1,\cdots,d\}$ and $\la>0$ we have
\beq
\label{1223}
&&\chi_{\la,i_1}^{(n)}(t,\bx)\\
&=&\sum\limits_{\bbi^{|2},
{\cal
    F}\in\cG_s,\bbsig}\int\cdots\int
    \vphi^{(n)}_{\bbi,\bbsig}(\bk_1,\cdots,\bk_n)P_{n-1}({\cal
    F})
\frac{Q({\cal
    F})}{\la+\sum\limits_{m\in A_n({\cal
    F})}|\bk_m|^{2\bt}}
    \prod\limits_{m\in
    A_n({\cal
    F})}\widehat{V}_{\si_m,i_m}(t,\bx,d\bk_m)\nonumber
\eeq
where  $\vphi^{(n)}_{\bbi,\bbsig}$ are some functions with $\sup{|\vphi^{(n)}_{\bbi,\bbsig}|}
\leq 1$
and
\beq
P_{n-1}({\cal
    F})&=&\prod\limits_{j=1}^{n-1}\frac{\sum\limits_{m\in A_j({\cal
    F})}|\bk_{m}|}{\la+\sum\limits_{m\in A_j({\cal
    F})}|\bk_m|^{2\bt}}
    \label{a.4}\\
\label{12203}
Q({\cal
    F})&=&\mathop{\prod\limits_{\widehat{mm'}\in E({\cal
    F})}}\bE\left[\widehat{V}_{\si_{m},i_{m}}(0,d\bk_{m})\widehat{V}_{\si_{m'},i_{m'}}(0,d\bk_{m'})\right]
\eeq
where $\bbsig$ is summed over all proper functions defined
over $\{1,\cdots,n\}$.
\end{lemma}

{\bf Remark.} Indexing terms of the sum in (\ref{1223})
with proper functions keeps track of differentiations
we have to perform in the process of constructing consequtive $\chi^{(n)}$-s.
The idea is as follows. The
 spectral  representation of the $n$-th order corrector contains terms which are multiple
 spectral
integrals with respect to the product measures of the form
\be
\label{01221}
\prod\limits_{m\in A_n({\cal F})}\widehat{V}_{\si_m,i_m}(t,\bx,d\bk_m).
\ee
Sigmas stand for the fact that there are two types of measures $\widehat{V}_0$, $\widehat{V}_1$
involved which are orthogonal to each other. In order to construct the next generation
corrector we have to apply nabla operator in $\bx$ to (\ref{01221}) using the Leibnitz rule.
This procedure consists in picking some factor $\widehat{V}_{\si_p,i_p}(t,\bx,d\bk_p)$ for
a certain $p\leq n$, belonging to $A_n({\cal F})$ and changing that term to
$\widehat{V}_{1-\si_p,i_p}(t,\bx,d\bk_p)$. The factor $\bk_p$ arising in the process
will be absorbed into the remaining expression. Finally to form $\chi^{(n+1)}$
we have to multiply yet the entire expression by
$\widehat{V}_{0,i_{n+1}}(t,\bx,d\bk_{n+1})$.
The manipulations we have just described are reflected by the rule of forming the proper
function of the consequtive generation described by formula (\ref{299}).

The utility of the correctors and convectors  depends on the
following lemma.
\begin{lemma}
\label{lemma1}
We have the recursive scheme, for $n\geq 1$,
\be
\label{124}
\int\limits_0^t\bU^{(n-1)}_{\ep}(s)ds=
-\chi_{\ep}^{(n)}(t)+\chi_{\ep}^{(n)}(0)+\int\limits_0^t\chi_{\ep}^{(n)}(s)ds+
\int\limits_0^t\bU^{(n)}_{\ep}(s)ds+
M_\ep^{(n-1)}(t).
\ee
where $M_\ep^{(0)}(t)$ is a Brownian Motion
with the covariance matrix
\be
\label{061104}
\bD_\ep=\int\limits_{R^d}\frac{|\bk|^{2\bt}a(|\bk|)}{|\bk|^{2\al+d-2}(|\bk|^{2\bt}+\ep^2)^2}\left(\bI-\frac{\bk\otimes\bk}{|\bk|^2}\right)d\bk.
\ee
and $M_\ep^{(n)}(t), n\geq1,$ are continuous  martingales with respect to
the filtration ${\cal
  V}_{-\infty,\frac{t}{\ep^2}}$.

For $\alpha+\beta<1$
  there exists $\ga(n)>0$ such that
\beq
\label{08121}
&&\lim\limits_{\ep\da} \ep^{-\ga(n)}\bE|\chi_{\ep}^{(n)}(t)|^2=0,\quad\mbox{for }n\geq1\\
&&\label{08122}
\lim\limits_{\ep\da}  \ep^{-\ga(n)}\bE|\bU^{(n)}_{\ep}(t)|^2=0,\quad
\hbox{for}\,\,n\geq 2[\bt]+1,~~[\beta]:
\mbox{the integer part of $\beta$}\\
\label{08123}
&&\lim\limits_{\ep\da} \bE|\ep^{-\ga(n)} M_\ep^{(n)}(t)|^p=0,\mbox{ for }n\geq 1
\eeq
\end{lemma}


We apply now Lemma~3 to prove Theorem \ref{theorem1}. By (\ref{124})  we have
\be
\label{90.10}
\bx_\ep(t)=M_\ep^{(0)}(t)+R_\ep(t)
\ee
with the remainder $R_\ep(t)$
\be
\label{30'}
R_\ep(t)=R_{s,\ep}(0)-R_{s,\ep}(t)+R_{m,\ep}(t)+R_{a,\ep}(t)
\ee
where
\[
R_{s,\ep}(t):=\sum\limits_{n=1}^{2[\beta]+1}\chi_\ep^{(n)}(t)
\]
is stationary,
\be
\nonumber
R_{a,\ep}(t):=\sum\limits_{n=1}^{2[\beta]+1}\int\limits_0^t\chi_{\ep}^{(n)}(s)ds+
\int\limits_0^t\bU^{(2[\beta]+1)}_{\ep}(s)ds
\ee
additive
and
\[
R_{m,\ep}(t):=\sum\limits_{n=1}^{2[\beta]}M_\ep^{(n)}(t)
\]
 a vanishing martingale.

Clearly, $M_\ep^{(0)}(t)$ converge,
 as continuous processes, to the Brownian motion with the covariance
 $\bD_0=\lim_{\ep\to 0}\bD_\ep$ (see, e.g., \cite{helland}).

To finish the proof it remains to show  tightness
$\lim\limits_{\ep\da}P[\sup\limits_{0\leq t\leq T}|R_\ep(t)|\geq
\delta]=0,\quad\forall T>0,\,\, \delta>0.
$
By Lemma \ref{lemma1} it suffices to show
$\lim\limits_{\ep\da}P[\sup\limits_{0\leq t\leq T}|R_{s,\ep}(t)|\geq
\delta]=0.
$
Since $R_{s,\ep}$ is stationary,
\be
\label{061107}
P[\sup\limits_{0\leq t\leq T}|R_{s,\ep}(t)|\geq
\delta]\leq\frac{T}{\ep^2}P[\sup\limits_{0\leq t\leq T\ep^2}|R_{s,\ep}(t)|\geq
\delta].
\ee
Using (\ref{90.10}), (\ref{30'}) we bound the right
side  of (\ref{061107})  by
\be
\label{061110}
\frac{T}{\ep^2}P[\sup\limits_{0\leq t\leq T\ep^2}[|R_{a,\ep}(t)|+|\bx_\ep(t)|]\geq
\delta]+\frac{T}{\ep^2}P[\sup\limits_{0\leq t\leq T\ep^2}[|R_{m,\ep}(t)|+|M^{(0)}_\ep(t)|]\geq
\delta]+
\ee
\[
\frac{T}{\ep^2}P[|R_{s,\ep}(0)|\geq
\delta].
\]
The last term in (\ref{061110})  is smaller than
\be
\label{061112}
\frac{\bE|R_{s,\ep}(0)|^pT}{\ep^2\delta^p}
\ee
by Chebychev's inequality.

Let $p$ be such that $\ga p>4$ where $\ga=\min\limits_{1\leq n\leq
  N}\ga(n)$. Expression (\ref{061112}) tends to $0$ as
$\ep\da$ by (\ref{08121}).

To prove that the second term of  (\ref{061110}) vanishes as $\ep\da$  we only need to show that for any $n\geq1$
\be
\label{08142}
\frac{T}{\ep^2}P[\sup\limits_{0\leq t\leq T\ep^2}|M^{(n)}_\ep(t)|]\geq
\delta]
\ee
tends to $0$ as $\ep\da$. Using (\ref{08122})  and Kolmogorov's inequality we
bound (\ref{08142}) by
\[
\frac{\bE|M^{(n)}_\ep(T)|^pT}{\ep^2\delta^p}\leq C\frac{\ep^{p\gamma(n)-2}}{T\delta^p}
\]
which vanishes with $\ep$.

The first term of (\ref{061110}) can be bounded by
\be
\label{061111}
C\left[\frac{\ep^2T}{\delta^2}\sum\limits_{n=1}^{2[\beta]+1}\bE|\chi_{\ep}^{(n)}(0)|^2\right.+
\frac{\ep^2T}{\delta^2}\bE|\ep^{2[\beta]}\bU^{(2[\beta]+1)}_{\ep^2}(0,\bze)|^2+
\left.\frac{T}{\ep^2}P[\int\limits_0^T\left|\vV(s,\ep\bx(s))\right|ds\geq\frac{\delta}{\ep}]\right].
\ee
Expression (\ref{061111}) tends to zero by Lemma
\ref{lemma1} and that
\beq
\limsup\limits_{\ep\da}\bE\left[\int\limits_0^T\left|\vV(s,\ep\bx(s))\right|ds\right]^p&\leq&
T^{p-1}\int\limits_0^T\bE\left|\vV(s,\ep\bx(s))\right|^pds\nonumber\\
&=&\bE\left|\vV(0,\bze)\right|^pT^{p}<+\infty. \nonumber
\eeq
For $p>2$ the last term in (\ref{061111})
vanishes with $\ep$.  Thus  the proof of Theorem \ref{theorem1}
is complete.

\section{Proof of Lemmas~2 and~3}

{\bf Proof of Lemma \ref{lemma2b}.} We prove  (\ref{1223}) by induction.
The case $n=1$  is obvious. We choose
$\vphi^{(0)}_{i_1}\equiv 1$. Suppose that the
result holds for $n$. Thus by (\ref{123}), (\ref{122}) and (\ref{a.0})
\beq
\nonumber
&&\chi_{\la,i_1}^{(n+1)}(t,\bx)\\
&&=\int\limits_t^{+\infty}e^{-\la(s-t)}\bE\left\{\vV(s,\bx)\cdot\nabla\chi_{\la,i_1}^{(n)}(s,\bx)\left|\right.
{\cal V}_{-\infty,t}\right\}ds\nonumber\\
\nonumber
&&=
\sum\limits_{\bbi^{|2},
{\cal
    F},\bbsig}\bE\left\{\int\limits_t^{+\infty}\int\cdots\int
    e^{-\la(s-t)}\vphi_{\bbi,\bbsig}^{(n)}(\bk_11,\cdots,\bk_n)P_{n-1}({\cal
    F})Q({\cal
    F})\frac{1}{\la+\sum\limits_{m\in A_n({\cal
    F})}|\bk_m|^{2\bt}}\right.\\
&&\left.\widehat{\vV}_0(s,\bx,d\bk_{n+1})\cdot\nabla\left\{
\prod\limits_{m\in
A_n({\cal
F})}\widehat{V}_{\si_{m},i_m}(s,\bx,d\bk_m)\right\}ds|
{\cal V}_{-\infty,t}\right\}.
\label{1224}
\eeq

To calculate (\ref{1224}) we decompose each $\widehat{V}_{\si,i}(s,\bx,d\bk)$
as $\widehat{V}_{\si,i}(s,\bx,d\bk)=\widehat{V}^0_{\si,i}(s,\bx,d\bk)
+\widehat{V}^1_{\si,i}(s,\bx,d\bk)$ where
\be
\label{a.3}
\widehat{V}^0_{\si,i}(s,\bx,d\bk)= e^{-|\bk|^{2\bt}(s-t)}\widehat{V}_{\si,i}
(t,\bx,d\bk)
\ee
is
the orthogonal projection onto
$L^2_{-\infty,t}$ and
\be
\label{a.6}
\widehat{V}^1_{\si,i}(s,\bx,d\bk)=(1-e^{-|\bk|^{2\bt}(s-t)})
\widehat{V}_{\si,i} (t,\bx,d\bk)
\ee
the
orthogonal complement.
We have
$\widehat{\vV}_{\si}(s,\bx,d\bk)=\widehat{\vV}^0_{\si}(s,\bx,d\bk)
+\widehat{\vV}^1_{\si}(s,\bx,d\bk)$.

Expression (\ref{1224}) becomes
\be
\label{1225}
\sum\limits_{\bbi^{|2},
{\cal
    F},\bbsig}\int\limits_t^{+\infty}\bE\left[\int\cdots\int
    e^{-\la(s-t)}\vphi_{\bbi,\bbsig} ^{(n)}(\bk_1,\cdots,\bk_n)
P_{n-1}({\cal
    F})Q({\cal
    F}){\cal K}({\cal
    F}) \left|\right.{\cal V}_{-\infty,t}\right]ds
\ee
with
\[
{\cal K}({\cal
    F}) := \frac{1}{\la+\sum\limits_{m\in A_n({\cal
    F})}|\bk_m|^{2\bt}}
\sum\limits_{
\varrho=\{\varrho_{j}\}\atop j\in A_n(\cF)\cup\{n+1\}}
\hspace*{.005in}~\widehat{\vV}_0^{\varrho_{n+1}}(s,\bx,d\bk_{n+1})\cdot\nabla\left\{\prod\limits_{m\in
    A_n({\cal
    F})}\widehat{V}_{\si_m,i_m}^{\varrho_{m}}(s,\bx,d\bk_m)\right\}.
    \]
The term corresponding to $\varrho_j\equiv 1$
vanishes,
as is shown in the following calculation,
\beq
\label{0311}
&&\bE\left\{\int\limits_t^{+\infty}\left[\int\cdots\int e^{-\la(s-t)}\vphi_{\bbi,\bbsig} ^{(n)}(\bk_1,\cdots,\bk_n)P_{n-1}({\cal
    F})Q({\cal
    F})\left(\la+\sum\limits_{m\in A_n({\cal
    F})}|\bk_m|^{2\bt}\right)^{-1}\right.\right.\\
&&\left.\left.\widehat{\vV}^{1}_0(s,\bx,d\bk_{n+1})
\cdot\nabla\left(\prod\limits_{m\in
    A_n({\cal
    F})}\widehat{V}_{\si_m,i_m}^{1}(s,\bx,d\bk_m)\right)\right]ds\right\}
    \nonumber\\
   & =&
\nabla\cdot\left\{\bE\int\limits_t^{+\infty}\left[\int\cdots\int
e^{-\la(s-t)}\vphi_{\bbi,\bbsig} ^{(n)}(\bk_1,\cdots,\bk_n)P_{n-1}({\cal
  F})Q({\cal
  F})\right.\right.\nonumber\\
&&\left.\left.\left(\la+\sum\limits_{m\in A_n({\cal
    F})}|\bk_m|^{2\bt}\right)^{-1}
    \widehat{\vV}^{1}_0(s,\bx,d\bk_{n+1})\prod\limits_{m\in
    A_n({\cal
    F})}\widehat{V}_{\si_m,i_m}^{1}(s,\bx,d\bk_m)\right]ds\right\}\nonumber\\
    &=&0\nonumber
    \eeq
    by homogeneity.

We use eqs. (\ref{a.2}) and (\ref{a.1})  to compute
\beq
&&\widehat{\vV}_0(s,\bx,d\bk_{n+1})\cdot\nabla\left\{\prod\limits_{m\in
A_n({\cal
F})}\widehat{V}_{\si_m,i_m}(s,\bx,d\bk_m)\right\}\label{a.7}\\
&=&\sum\limits_{m'\in
A_n({\cal
F})}\bk_{m'}\cdot\widehat{\vV}_{\si_{n+1}^{m'}}(s,\bx
,d\bk_{n+1})\prod\limits_{m\in
A_n({\cal
F})}\widehat{V}_{\si_m^{m'},i_m}(s,\bx,d\bk_m)
\nonumber
\eeq
where $\bbsig^{m'}$ is a proper function obtain from $\bbsig$ via
procedure described by (\ref{299}) with $p=m'$.
\commentout{\be
\label{46'}
\si_{m}^{m'}:=\left\{
\begin{array}{ll}
1-\si_{m'} &\mbox{if\quad}m'=m\\
\si_m &\mbox{otherwise.}
\end{array}\right.
\ee}

By (\ref{a.3}), (\ref{a.6}), (\ref{a.7}) and the definition
(\ref{a.4}), (\ref{1225}) further reduces to
\be
\label{061605}
\sum\limits_{\bbi^{|2},
{\cal
F},\bbsig}
\int\limits_t^{+\infty}\left\{\int\cdots\int\sum\limits_{i_{n+1}=1}^d
\sum\limits_{m'\in A_n({\cal F})}\sum\limits_{{\cal F}'}
\vphi_{\bbi,\bbsig} ^{(n)}
\frac{k_{m', i_{n+1}}}{\sum\limits_{m\in A_n({\cal F})}
|\bk_{m}|}\right.
\ee
\[
\exp\{-(\la+\sum\limits_{m\in A({\cal
 F}')}|\bk_m|^{2\bt})\left(s-t\right)\}P_{n}({\cal
 F})Q({\cal
 F})
 \prod\limits_{m\in
  A({\cal F}')}\widehat{V}_{\si_m^{m'},i_m}(t,\bx,d\bk_m)
  \]
  \[
\left.\mathop{\prod\limits_{\widehat{pq}\in E({\cal
F}')}}\left[1-e^{-(|\bk_{p}|^{2\bt}+|\bk_{q}|^{2\bt})(s-t)}\right]\bE\left[
\widehat{V}_{\si_{p}^{m'},i_{p}}
(0,\bze,d\bk_{p})\widehat{V}_{\si_{q}^{m'
 },i_{q}}(0,\bze,d\bk_{q})\right]\right\}.
 \]
Here $\sum\limits_{{\cal F}'}$
  denotes the summation over all
  incomplete Feynman diagrams ${\cal F}'$ based on the
  set $A_n({\cal F})\cup \{n+1\}$.

The proof of the lemma is complete with
\[
\vphi_{\bbi,\bbsig^{m'}}^{(n+1)}(\bk_1,\cdots,\bk_{n+1}):=
\vphi_{\bbi,\bbsig} ^{(n)}(\bk_1,\cdots,\bk_{n})
\frac{k_{m',i_{n+1}}}{\sum\limits_{m'\in A_n({\cal F})}
  |\bk_{m'}|}\cdot
\]
\[
\left(\la+\sum\limits_{m\in A({\cal F}')}|\bk_m|^{2\bt}\right)\int\limits_t^{+\infty}
e^{-(\la+\sum\limits_{m\in A({\cal F}')}|\bk_m|^{2\bt})\left(s-t\right)}
\mathop{\prod\limits_{\widehat{pq}\in E({\cal
    F}')}}\left[1-e^{-(|\bk_{p}|^{2\bt}+|\bk_{q}|^{2\bt})(s-t)}\right]ds\mbox{\kropa}
\]

{\bf Proof of Lemma \ref{lemma1}}.
\commentout{
Throughout this proof we suppress the subscript $\ep$ and write
$\bU^{(n)}$ and  $\chi^{(n)}$.}

We recall the {\em pseudogenerator} of the path
corrector
\[
\cL\chi^{(n)}_{\ep}(t):=\lim\limits_{\delta\da}
\frac{\bE\left[\chi^{(n)}_{\ep}(t+\delta)\left|\right.{\cal
  V}_{-\infty,\frac{t}{\ep^2}}\right]-\chi^{(n)}_{\ep}(t)}{\delta}
\]
in the $L^1$ sense.
By a finite Taylor expansion
we have
\be
\label{01205}
\chi^{(n)}_{\ep}(t+\delta)=
\ep^n\chi_{\ep^2}^{(n)}(\frac{t+\delta}{\ep^2},\bx_\ep(t))+R_\delta(t),
\ee
with
\[
R_\delta(t)=
\ep^{n-1}\bU^{(n)}_{\ep^2}
(\frac{t}{\ep^2},\bx_{\ep}(t))\delta+o(\delta)
\]
where $o(\delta)$ is in the
$L^p$ sense.

>From (\ref{123}) and
(\ref{122}) we have
\beq
\ep^n\bE\left[\chi_{\ep^2}^{(n)}(\frac{t+\delta}{\ep^2},\bx_\ep(t))\left|\right.{\cal
 V}_{-\infty,\frac{t}{\ep^2}}\right]&=&
\ep^n\int\limits_{\frac{t+\delta}{\ep^2}}^{+\infty}
e^{-\ep^2(s-\frac{t+\delta}{\ep^2})}\bE\left[\bU^{(n-1)}_{\ep^2}
(s,\bx_\ep(t))\left|\right.{\cal
 V}_{-\infty,\frac{t}{\ep^2}}\right]
\label{a.5}\\
&=&\chi^{(n)}(t)+
\delta\left[\chi^{(n)}(t)-\ep^{n-2}\bU^{(n-1)}_{\ep^2}
(\frac{t}{\ep^2},\bx_\ep(t))\right]+o(\delta).\nonumber
\eeq

Thus,
\be
\label{01203}
\cL\chi^{(n)}_{\ep}(t)= \chi^{(n)}_{\ep}(t)-\bU^{(n-1)}_{\ep}(t)
+\bU^{(n)}_{\ep}(t)=\cL_{\ep}(t;\bx_{\ep}(t))+\bU^{(n)}_{\ep}(t).
\ee
Here $\cL_{\ep}(t;\bx)$, $t\geq0$ denotes the pseudogenerator of the process
$\ep^n\chi^{(n)}(\frac{t}{\ep^2},\bx)$, $t\geq0$ calculated with respect
to the filtration ${\cal V}_{-\infty,\frac{t}{\ep^2}}$, $t\geq0$
for a fixed $\bx\in R^d$. The processes
\be
\label{01206}
M_{\ep}^{(n)}(t;\bx):=\ep^n\chi^{(n)}(\frac{t}{\ep^2},\bx)-
\int\limits_0^t\cL_{\ep}(s;\bx)ds,
~\mbox{ for a fixed }\bx\in R^d
\ee
and
\be
\label{3111}
M^{(n)}_{\ep}(t)=\chi^{(n)}_{\ep}(t)-
\int\limits_0^t \cL\chi^{(n)}_{\ep}(s)ds
\ee
are both martingales (cf. Section \ref{sec2}).

The quadratic variation of the first martingale
can be calculated with the help of Proposition \ref{prop5} and it is given by
\[
<M_{\ep,i}^{(n)}(\cdot;\bx)M_{\ep,j}^{(n)}(\cdot;\bx)>_t=
\int\limits_0^t W_{\ep,i,j}^{(n)}(s,\bx)ds,
\]
where
\be
\label{012016}
W_{\ep,i,j}^{(n)}(s,\bx)=
\ee
\[
\ep^{2n-2}\sum~'
\commentout{\mathop{\sum\limits_{\bbi_{|1}=i,
{\cal
    F},\bbsig}\sum\limits_{\bbi'_{|1}=j,
    {\cal
        F}',\bbsig'}}\mathop{\sum\limits_{q\in A_n({\cal F})}}
\sum\limits_{q'\in A_n({\cal F}') }}
\int\cdots\int
    \vphi^{(n)}_{\bi,\bbsig}\vphi^{(n)}_{\bi',\bbsig'}P_{n-1}({\cal
    F})Q({\cal
    F})P_{n-1}({\cal
    F}')Q({\cal
    F}')
\]
\[
\frac{a(|\bk_q|)|\bk_q|^{2\bt}}{|\bk_q|^{2\al+d-2}(\ep^2+\sum\limits_{m\in A_n({\cal
    F})}|\bk_m|^{2\bt})(\ep^2+\sum\limits_{m\in A_n({\cal
    F}')}|\bk_m'|^{2\bt})}\left(\delta_{i_q,i_{q'}'}-\frac{k_{q,i_q}k_{q',i_{q'}'}}
{|\bk_q|^2}\right)
\]
\[
\delta_{\si_q,\si_{q'}'}\delta(\bk_q-\bk_{q'})d\bk_qd\bk_{q'}\prod\limits_{m\in
    A_n({\cal
    F}),m\not=q}\widehat{V}_{\si_{m},i_m}(\frac{t}{\ep^2},\bx,d\bk_m)\prod\limits_{m\in
    A_n({\cal
    F}'),m\not=q'}\widehat{V}_{\si_{m}',i_m'}
    (\frac{t}{\ep^2},\bx,d\bk_m').
\]
The summation on the right hand side of the above expression extends over all multiindices of
length $n$,
$\bbi,\bbi'$ whose first components are $i$ and $j$ correspondingly, all pairs of 
diagrams ${\cal F}$, ${\cal F}'\in {\cal G}_s$ made of vertices of ${\cal S}=\{1,\cdots,n\}$,
all pairs of proper functions $\bbsig,\bbsig':{\cal S}\rightarrow\{0,1\}$ and all pairs of $q,q'$
 belonging to $A_n({\cal F})$, $A_n({\cal F}')$ respectively.
We have denoted here $\bk_m=(k_{m,1},\cdots,k_{m,d})$.

As for the martingale given by (\ref{3111}), after
elementary calculations we get
\be
<M^{(n)}_{\ep,i},M^{(n)}_{\ep,j}>_t=\int\limits_0^t
\cL\left\{M^{(n)}_{\ep,i}M^{(n)}_{\ep,j}\right\}(s)ds
\label{061601}
\ee
with
\be
\label{012010}
\cL\left\{M^{(n)}_{\ep,i}M^{(n)}_{\ep,j}\right\}(s)=\cL\left\{\chi^{(n)}_{\ep,i}
\chi^{(n)}_{\ep,j}\right\}(s)
-\cL\chi^{(n)}_{\ep,i}(s)\chi^{(n)}_{\ep,j}(s)-
\chi^{(n)}_{\ep,i}(s)\cL\chi^{(n)}_{\ep,j}(s).
\ee
\commentout{
\ep^{2n}\lim\limits_{\delta\da}\frac{\bE\left[\chi^{(n)}_i(\frac{t+\delta}{\ep^2},\bx_\ep(t))\chi^{(n)}_j(\frac{t+\delta}{\ep^2},\bx_\ep(t))\left|\right.{\cal
V}_{-\infty,\frac{t}{\ep^2}}\right]-
\chi^{(n)}_i(\frac{t}{\ep^2},\bx_\ep(t))\chi^{(n)}_j(\frac{t}{\ep^2}
,\bx_\ep(t))}{\delta}}
Calculating precisely as in (\ref{01205})
we obtain that the first term on the right hand side of (\ref{012010}) equals
\[
\cL_{\ep,i,j}'(s;\bx_\ep(s))+
\chi^{(n)}_{\ep,i}(s)U^{(n)}_{\ep,j}(s)+U^{(n)}_{\ep,i}(s)\chi^{(n)}_{\ep,j}(s).
\]
Here $\cL_{\ep,i,j}'(s;\bx)$, $s\geq0$ denotes the pseudogenerator of
the process $\ep^{2n}\chi^{(n)}_i(\frac{s}{\ep^2},\bx)
\chi^{(n)}_j(\frac{s}{\ep^2},\bx)$, $s\geq0$ with $\bx\in R^d$ fixed.
Taking  into account (\ref{01203}) we obtain that
\be
\label{012015}
\cL\left\{M^{(n)}_{\ep,i}M^{(n)}_{\ep,j}\right\}(s)=
W_{\ep,i,j}^{(n)}(s,\bx_\ep(s)).
\ee

In particular,  the above calculation shows that
$M^{(0)}$ is a Brownian Motion whose covariance matrix equals (\ref{061104}).

\subsection{Proof of eqs. (\ref{08121}) - (\ref{08123}).}
We break the proof into two cases: $\al+2\bt>2$ and $\al+2\bt\leq 2$.

{\bf The case of $\al+2\bt>2$}.
>From $\al+\bt<1$  and $\al+2\bt>2$ it follows that $\bt>1$.
For $2\bt>1$,
there exists a constant $C(n,\bt)$, depending only on
$n$ and $\bt$, such that, for any $m_j\in A_j({\cal
    F})$,
\be
\label{012011}
\frac{\sum\limits_{m\in A_j({\cal
    F})}|\bk_m|}{\ep^2+\sum\limits_{m\in A_j({\cal
    F})}|\bk_m|^{2\bt}}\leq
    C\frac{|\bk_{m_0}|+\ep^{\frac{1}{\bt}}}{\ep^2+|\bk_{m_0}|^{2\bt}}.
\ee
So,
\be
\label{224}
P_{n-1}({\cal F})\leq C\prod\limits_{j=1}^{n-1}
\frac{|\bk_{m_j}|+\ep^{\frac{1}{\bt}}}{\ep^2+|\bk_{m_j}|^{2\bt}}
\ee
for any $m_j\in A_j({\cal
F})$. Let
$m_j:=j$ if $j$ is {\em not} the right endpoint of an edge.
Otherwise, let $m_j$ be the closest free vertex to the left
of the edge whose right endpoint is $j$.

>From (\ref{224}) we have
\be
\label{1226}
\bE\left[\chi_{i_1}^{(n)}(0)\right]^2\leq
\ee
\[
C\ep^{2n-2/\beta}
\sum\limits_{\bi_{|2}}\sum\limits_{{\cal F}\in\cG_s,\si_l}\int\cdots\int
\int\cdots\int\prod\limits_{m\in A_n({\cal F})}\left[\frac{(|\bk_{m}|+
\ep^{\frac{1}{\bt}})(|\bk_{m}'|+\ep^{\frac{1}{\bt}})}{\left(\ep^2+
|\bk_{m}|^{2\bt}\right)\left(\ep^2+|\bk_{m}'|^{2\bt}\right)}\right]^{1+\al_m}
\]
\[
\mathop{\prod\limits_{\widehat{pq}\in E({\cal F})}}\frac{(|\bk_{p}|+
\ep^{\frac{1}{\bt}})(|\bk_{p}'|+\ep^{\frac{1}{\bt}})a(|\bk_{p}|)
a(|\bk_{p}'|)d\bk_{p}d\bk_{p}'}{\left(\ep^2+|\bk_{p}|^{2\bt}\right)
\left(\ep^2+|\bk_{p}'|^{2\bt}\right)(|\bk_{p}||\bk_{p}'|)^{d+2\al-2}}
\]
\[
\left|\bE\left[\prod\limits_{m\in A_n({\cal F})}\widehat{V}_{i_{m}}(0,d\bk_{m})
\widehat{V}_{i_{m}}(0,d\bk_{m}')\right]\right|,
\]
where $\al_m\geq0$ is the number of
{\em left } endpoints between a free vertex $m\in A_n(\cF)$
and the next free vertex $m'\in A_n({\cal F})$.
Here and below the
variables $\bk_p'$-s denote the second set of wavenumbers arising as a
result of squaring of the right side of the expansion formula for $\chi^{(n)}$
(\ref{1223}).
We also use
 the notation $\bk_{n+p}=\bk_p'$, thinking that
 the second copy of the diagram is based on the shifted set
 $\{n+1,n+2,...,2n\}$.

Adopting the notation, $k_p=|\bk_p|$, subsequently and
computing the expectation on the right side of (\ref{1226}) we get that
\be
\label{301}
\bE\left[\chi_{i_1}^{(n)}(0)\right]^2\leq
C\ep^{2n-2/\bt}\sum\limits_{{\cal
    F}\in\cG_s}\sum\limits_{{\cal F}'}
\int\limits_0^K\cdots\int\limits_0^K
\prod\limits_{m\in A_n({\cal
F})}\left[\frac{(k_{m}+\ep^{\frac{1}{\bt}})(k_{m}'+\ep^{\frac{1}{\bt}})}{\left(\ep^2+k_{m}^{2\bt}\right)\left(\ep^2+(k_{m}')^{2\bt}\right)}\right]^{1+\al_m}
\ee
\[
\mathop{\prod\limits_{\widehat{pq}\in E({\cal
    F}')}}\frac{\delta(k_{p}-k_{q})dk_{p}dk_{q}}{k_{p}^{2\al-1}}
\mathop{\prod\limits_{\widehat{pq}\in E({\cal F})}}\frac{(k_{p}+\ep^{\frac{1}{\bt}})(k_{p}'+\ep^{\frac{1}{\bt}})dk_{p}dk_{p}'}{\left(\ep^2+k_{p}^{2\bt}\right)\left(\ep^2+(k_{p}')^{2\bt}\right)(k_{p}k_{p}')^{2\al-1}}
\]
where
${\cal F}'$ is summed over $\cG_c(A_n(\cF))\cup\{n+A_n(\cF)\}$.

Note that
\be
\label{302}
c_n+2\sum\limits_{m\in A_n({\cal
F})}\al_m=n,
\ee
where $c_n$ denotes the cardinality of the set $A_n({\cal
F})$.

For $\alpha+2\beta>2$ the integrals with respect to the
variables $k_{p}$, $p$ - the vertices of the left edges of ${\cal F}'$, become
unbounded when $\ep\da$. Indeed as a result of integrating out the Dirac's delta
functions we obtain terms whose singularity at $0$ are becoming, when $\ep\da$
of the same order of magnitude as $\frac{1}{k_p^{(2+\al_p+\al_q)(2\bt-1)+2\al-1}}>>
\frac{1}{k_p^{4\bt+2\al-3}}$. This fact gives rise to divergence
of the corresponding improper integral at $0$ in the above regime of parameters $\al,\bt$.
Since $\bt>1$, otherwise $\al+2\bt<2$ and we are not considering that case
here this singularity can b
Then changing variables by setting
$\tilde{k}_{p}=k_{p}\ep^{-\frac{1}{\bt}}$ we obtain using (\ref{302}) that
\be
\label{2505}
\bE\left[\chi_{i_1}^{(n)}(0)\right]^2\leq
C\sum\limits_{{\cal
F}}\sum\limits_{{\cal F}'}\ep^{\frac{n-2+c_n+2c_n(1-\al-\bt)}{\bt}}
\int\limits_0^{K}\cdots\int\limits_0^{K}\mathop{\prod\limits_{\widehat{pq}
\in E({\cal F})}}\frac{(k_{p}+\ep^{\frac{1}{\bt}})(k_{p}'+
\ep^{\frac{1}{\bt}})dk_{p}dk_{p}'}{\left(\ep^2+k_{p}^{2\bt}\right)
\left(\ep^2+(k_{p}')^{2\bt}\right)(k_{p}k_{p}')^{2\al-1}}
\ee
\[
\int\limits_0^{K\ep^{-\frac{1}{\bt}}}\cdots\int\limits_0^{K\ep^{-\frac{1}{\bt}}}\mathop{\prod}\limits_{\widehat{pq}\in
E({\cal F}')}\frac{(\tk_{p}+1)^{2+\al_{p}+\al_{q}}d\tk_{p}}{\left(\tk_{p}^{2\bt}+1\right)^{2+\al_{p}+\al_{q}}\tk_{p}^{2\al-1}}.
\]
The first multiple integral appearing on the right hand side of (\ref{2505})
remains bounded when $\ep\da$ thanks to $\al +\bt<1$. The utmost right
hand side factor however becomes then a product of improper integrals.
The  rate of decay of the integrand corresponding
to index $p$ in that integral is the same as for
\[
\frac{1}{\tk_{p}^{(2\bt-1)(2+\al_p+\al_q)+2\al-1}}<<
\frac{1}{\tk_{p}^{4\bt+2\al-3}},\quad\mbox{for }\tk_{p}>>1.
\]
Since $\alpha+2\beta>2$
the  integrals  involved are convergent.
(\ref{08121}) follows then from the fact that $c_n\geq 1$.

{\bf  Proof of (\ref{08123})}. From Lemma \ref{lemma2b} and (\ref{122})
we have that
\be
\label{303}
U^{(n)}_{i_1}(t,\bx)=
\sum\limits_{\bi_{|1}=i_1}\sum\limits_{{\cal
    F },\bbsig}\int\cdots\int
    \vphi^{(n+1)}_{\bi,\bbsig}(\bk_1,\cdots,\bk_{n+1})P_{n}({\cal
    F})Q({\cal
    F})
\ee
\[
\prod\limits_{m\in
    A_{n}({\cal
    F})\cup\{n+1\}}\widehat{V}_{\si_{m},i_m}(t,\bx,d\bk_m),~~i_1=1,\cdots,d,
\]
where $\bU^{(n)}=(U^{(n)}_{1},\cdots,U^{(n)}_{d})$ and
  $|\vphi^{(n+1)}_{\bi,\bbsig}|\leq 1$.
  $\bE|\bU_n(0,\bze)|^2$ can be estimated in the same way as
  $\chi^{(n)}$. We outline the argument in the following.

First
(\ref{224}) and the choice of $m_j$
we get a bound, analogous to (\ref{1226}):
\[
\ep^{2(n-1)}\bE|\bU^{(n)}(0,\bze)|^2\leq
\]
\[
C\ep^{2(n-1)}
\sum\limits_{\bi}\sum\limits_{{\cal F}}\int\cdots\int\prod\limits_{m\in A_n({\cal F})\cup{\{n+1\}}}\left[\frac{(|\bk_{m}|+\ep^{\frac{1}{\bt}})(|\bk_{m}'|+\ep^{\frac{1}{\bt}})}{\left(\ep^2+|\bk_{m}|^{2\bt}\right)\left(\ep^2+|\bk_{m}'|^{2\bt}\right)}\right]^{1+\al_m}
\]
\[
\mathop{\prod\limits_{\widehat{pq}\in E({\cal F})}}\frac{(|\bk_{p}|+\ep^{\frac{1}{\bt}})(|\bk_{p}'|+\ep^{\frac{1}{\bt}})a(|\bk_{p}|)a(|\bk_{p}'|)d\bk_{p}d\bk_{p}'}{\left(\ep^2+|\bk_{p}|^{2\bt}\right)\left(\ep^2+|\bk_{p}'|^{2\bt}\right)(|\bk_{p}||\bk_{p}'|)^{d+2\al-2}}
\]
\[
\left|\bE\left[\prod\limits_{m\in A_n({\cal F})\cup{\{n+1\}}}\widehat{V}_{i_{m}}(0,d\bk_{m})\widehat{V}_{i_{m}}(0,d\bk_{m}')\right]\right|.
\]
Here we define $\al_{n+1}=-1$. Repeating the calculations
leading to (\ref{2505}) and noting
\be
\label{70}
c_n+2\sum\limits_{m\in A_n({\cal
 F})\cup{n+1}}\al_m=n-2,
\ee
we obtain using (\ref{70})
\be
\label{052501}
\ep^{2(n-1)}\bE|\bU^{(n)}(0,\bze)|^2\leq
\ee
\[
C\sum\limits_{{\cal
    F}\in\cG_s}\sum\limits_{{\cal F}'}\ep^{\frac{n+c_n+2(c_n+1)(1-\al-\bt)}{\bt}}
\int\limits_0^{K}\cdots\int\limits_0^{K}\mathop{\prod\limits_{\widehat{pq}
\in E({\cal F})}}\frac{(k_{p}+\ep^{\frac{1}{\bt}})(k_{p}'+\ep^{\frac{1}{\bt}})dk_{p}dk_{p}'}{\left(\ep^2+k_{p}^{2\bt}\right)\left(\ep^2+(k_{p}')^{2\bt}\right)(k_{p}k_{p}')^{2\al-1}}
\]
\[
\int\limits_0^{K\ep^{-\frac{1}{\bt}}}\cdots\int\limits_0^{K\ep^{-\frac{1}{\bt}}}\mathop{\prod\limits_{p <q }}\limits_{\widehat{pq}\in E({\cal F}')}\frac{(\tk_{p}+1)^{2+\al_{p}+\al_{q}}d\tk_{p}}{\left(\tk_{p}^{2\bt}+1\right)^{2+\al_{p}+\al_{q}}\tk_{p}^{2\al-1}}
\]
with ${\cal F}'\in \cG_c(A_n({\cal F})\cup\{n+1\}\cup
\{n+1+A_n({\cal F})\cup\{n+1\}\}$.

The  multiple
integral appearing in the last factor on the right hand side of
(\ref{052501}) becomes improper when we pass to the limit with
$\ep\da$. However the integrals involving variables with both $\al_{p},\al_{q}$ nonnegative
are convergent since $\al+2\bt >2$. The only possible divergence
comes therefore from the factors involving $\al_{n+1}=-1$. The worst
type of divergence could take place
for the exponent $\al_{n+1}$ appearing twice in such a factor.
This case corresponds  to an occurrence of the edge $(n+1,2n+2)$.
The end result is then the following asymptotic
\[
\ep^{2(n-1)}\bE|\bU^{(n)}(0,\bze)|^2\leq C\sum\limits_{{\cal
    F}}\ep^{\frac{n+c_{n}+2(c_{n}+1)(1-\al-\bt)-2+2\al}{\bt}}\leq C\ep^{\frac{n+3-2\al-4\bt}{\bt}}.
\]
The assertion of the lemma follows for
$n\geq 2[\bt]$.

{\bf Proof of (\ref{08123})}. We can use  (\ref{012015}) to represent the quadratic
variation of the martingale in question. After taking its
expectation and using \cite{port-stone} we get the following estimate, valid for any $p\geq1$
\be
\label{08141}
\bE|M_\ep^{(n)}(t)|^p\leq C t^{p}\sum\limits_{i=1}^d
\left(\bE W^{(n)}_{\ep,i,i}(0,\bze)\right)^{p/2}
\ee
where $W^{(n)}_{\ep,i,i}$ is defined by (\ref{012016}).
On the other hand using that formula one gets estimating in  the precisely  same way as it was done for (\ref{08121}) and (\ref{08122}) that
\be
\label{052502}
\bE W^{(n)}_{\ep,i,i}(0,\bze)^2\leq C\ep^{2n-2}\sum_{\bbi}
\sum_{{\cal
    F}\in \cG_s}\sum\limits_{m'\in
    A_n({\cal
    F})}\sum_{{\cal
    F}'}\mathop{\underbrace{\int\limits_{0}^K\cdots\int\limits_{0}^K}}\limits_{2n~~integrals}
\ee
\[
\mathop{\prod\limits_{m\in A_n({\cal
    F})}}\limits_{m\not=m'}\left[\frac{(k_{m}+\ep^{\frac{1}{\bt}})(k_{m}'+\ep^{\frac{1}{\bt}})}{\left(k_{m}^{2\bt}+\ep^2\right)\left((k_{m}')^{2\bt}+\ep^2\right)}\right]^{1+\al_m} \frac{k_{m'}^{2\bt}}{k_{m'}^{2\al-1}(k_{m'}^{2\bt}+\ep^2)^2}
\left(\frac{k_{m'}+\ep^{\frac{1}{\bt}}}{k_{m'}^{2\bt}+\ep^2}\right)^{2\al_{m'}}
\]
\[
\mathop{\prod\limits_{\widehat{pq}\in E({\cal F}')}}\frac{\delta(k_{p}-k_{q})dk_{p}dk_{q}}{k_{p}^{2\al-1}}
\mathop{\prod\limits_{\widehat{pq}\in E({\cal F})}}\frac{(k_{p}+\ep^{\frac{1}{\bt}})(k_{p}'+\ep^{\frac{1}{\bt}})dk_{p}dk_{p}'}{\left(\ep^2+k_{p}^{2\bt}\right)\left(\ep^2+(k_{p}')^{2\bt}\right)(k_{p}k_{p}')^{2\al-1}}
\]
with ${\cal
 F}'\in\cG_c(
 A_n({\cal
 F})\cup\{n+A_n(\cF)\}\
\{m', m'+n\}$.

Two situations arise: $\al_{m'}=0$ and $\al_{m'}\geq 1$.
For $\al_{m'}=0$, with an appropriate
change of variables we rewrite
the term in the above summations
\[
\ep^{\frac{2(c_n-1)(1-\al-\bt)+c_n+n-2}{\bt}}\int\limits_{0}^K\cdots\int\limits_{0}^K\mathop{\prod\limits_{\widehat{pq}\in E({\cal F})}}\frac{(k_{p}+\ep^{\frac{1}{\bt}})(k_{p}'+\ep^{\frac{1}{\bt}})dk_{p}dk_{p}'}{\left(\ep^2+k_{p}^{2\bt}\right)\left(\ep^2+(k_{p}')^{2\bt}\right)(k_{p}k_{p}')^{2\al-1}}
\]
\[
\int\limits_{0}^K\frac{k_{m'}^{2\bt}}{k_{m'}^{2\al-1}(k_{m'}^{2\bt}+\ep^2)^2}
dk_{m'}
\int\limits_{0}^{K\ep^{-\frac{1}{\bt}}}\cdots\int\limits_{0}^{K\ep^{-\frac{1}{\bt}}}\mathop{\prod\limits_{m\in A_n({\cal
    F})}}\limits_{m'\not =m}\left[\frac{(k_{m}+1)(k_{m}'+1)}{\left(k_{m}^{2\bt}+1\right)\left((k_{m}')^{2\bt}+1\right)}\right]^{1+\al_m}
\]
\[
\mathop{\prod\limits_{\widehat{pq}\in E({\cal F}')}}\frac{\delta(k_{p}-k_{q})dk_{p}dk_{q}}{k_{p}^{2\al-1}}\sim
C\ep^{\frac{2(c_n-1)(1-\al-\bt)+c_n+n-2}{\bt}}.
\]

For
$\al_{m'}\geq 1$, a similar calculation leads to
\[
\ep^{2n-2}\int\limits_{0}^K\cdots\int\limits_{0}^K\mathop{\prod\limits_{\widehat{pq}\in E({\cal F})}}\frac{(k_{p}+\ep^{\frac{1}{\bt}})(k_{p}'+\ep^{\frac{1}{\bt}})dk_{p}dk_{p}'}{\left(\ep^2+k_{p}^{2\bt}\right)\left(\ep^2+(k_{p}')^{2\bt}\right)(k_{p}k_{p}')^{2\al-1}}
\]
\[
\int\limits_{0}^{K\ep^{-\frac{1}{\bt}}}\cdots\int\limits_{0}^{K\ep^{-\frac{1}{\bt}}}\frac{k_{m'}^{2\bt}}{k_{m'}^{2\al-1}(k_{m'}^{2\bt}+\ep^2)^2}
\mathop{\prod\limits_{m\in A_n({\cal
    F})}}\limits_{m'\not =m}
\left[\frac{(k_{m}+1)(k_{m}'+1)}{\left(k_{m}^{2\bt}+1\right)\left((k_{m}')^{2\bt}+1\right)}\right]^{1+\al_m}
\]
\[
\mathop{\prod\limits_{\widehat{pq}\in E({\cal F}')}}\frac{\delta(k_{p}-k_{q})dk_{p}dk_{q}}{k_{p}^{2\al-1}}dk_{m'}\sim
C\ep^{\frac{2c_n(1-\al-\bt)+c_n+n-2}{\bt}}.
\]
In either situation, each term is of order $O(\ep^{\ga})$ for some $\ga>0$, and
$\forall n\geq 2$.

The special  case of $\alpha+2\beta=2$ can be similarly analyzed. This time,
some integrals would  diverge logarithmically, as $\vep\to 0$, but would
be controlled by factors of order $O(\ep^\gamma)$ for some $\gamma>0$.

{\bf The case of $\al+2\bt<2$.} We shall prove (\ref{08121})-(\ref{08123}) only
for $n=1$. This is sufficient for
the proof of Theorem \ref{theorem1}. The general case
requires evaluating multiple spectral integrals as in the previous case
of $\al+2\bt\geq 2$, only much more easily
because the corresponding improper
integrals are convergent.

We have
\[
\bE\left|\chi^{(1)}_\ep(t)\right|^2=\bE\left|\chi^{(1)}_\ep(0)\right|^2\leq\ep^2\int\limits_0^{+\infty}\frac{{\cal
    E}(k)}{(\ep^2+k^{2\bt})^2}dk\leq
\]
\[
C\ep^{\frac{2(1-\al-\bt)}{\bt}}\int\limits_0^{+\infty}\frac{~~dk}{(1+k^{2\bt})^2k^{2\al-1}}
\]
which vanishes, as $\ep\to 0$, for $\alpha+\beta<1$. Also,
\be
\label{221}
\bE\left|\bU^{(1)}_{\la}(0,\bze)\right|^2\leq
C\bE\left|\vV(0,\bze)\right|^2\bE\left|\nabla\chi^{(1)}_{\la}(0,\bze)\right|^2
\ee
and, consequently,
\be
\label{222}
\limsup\limits_{\la\da}\bE\left|\bU^{(1)}_{\la}(0,\bze)\right|^2<+\infty
\ee
as
\be
\label{40'}
\bE\left|\nabla\chi^{(1)}_{\la}(0,\bze)\right|^2\leq
\int\limits_{R^d}\frac{|\bk|^2\cE(|\bk|)}{(\lambda+|\bk|^{2\bt})^2}d\bk\leq
C(1+ \lambda^{2-\alpha-2\beta\over \beta})
\ee
is bounded, as $\lambda\to 0,$ for $\alpha+2\beta<2$. The
inequality (\ref{222}) holds for any $p$-th moment of $\bU_{\la}^{(1)}$, with
$p>0$ because of Gaussianity of the velocity field.

We also have
\[
\bE\left|\int\limits_0^t\bU^{(1)}_{\ep^2}(\frac{s}{\ep^{2}},\bbxo)ds\right|^2=
2\int\limits_0^tds\int\limits_0^s\bE\left\{\bU^{(1)}_{\ep^2}(\frac{s}{\ep^{2}},\bbxo)\cdot\bU^{(1)}_{\ep^2}(\frac{s_1}{\ep^{2}},\bmx_\ep(s_1))\right\}ds_1
\]
\be
\label{31.12}
=2\int\limits_0^tds\int\limits_0^s\bE\left\{\bU^{(1)}_{\ep^2}(\frac{s}{\ep^{2}},\bmx_\ep(s_1))\cdot\bU^{(1)}_{\ep^2}(\frac{s_1}{\ep^{2}},\bmx_\ep(s_1))\right\}ds_1+
\ee
\[
\frac{2}{\ep}\int\limits_0^tds\int\limits_0^sds_1\int\limits_{s_1}^s\bE\left\{\left(\vV(\frac{{s_2}}{\ep^{2}},\bx_\ep({s_2}))\cdot\nabla\bU^{(1)}_{\ep^2}(\frac{s}{\ep^{2}},\bx_\ep({s_2}))\right)\cdot\bU^{(1)}_{\ep^2}(\frac{{s_1}}{\ep^{2}},\bmx_\ep(s_1))\right\}d{s_2}
\]

Because of stationarity  of
$\bU^{(1)}_{\ep^2}$ along the trajectory (see \cite{port-stone})
the first term of (\ref{31.12}) equals
\[
2\int\limits_0^tds\int\limits_0^s\bE\left\{\bU^{(1)}_{\ep^2}(\frac{{s_1}}{\ep^{2}},\bze)\cdot\bU^{(1)}_{\ep^2}(0,\bze)\right\}d{s_1}=
\]
\be
\label{1.5}
2\int\limits_0^tds\int\limits_0^s\bE\left\{\left(\vV(\frac{{s_1}}{\ep^{2}},\bze;0)\cdot\nabla\chi^{(1)}_{\ep^2}(\frac{{s_1}}{\ep^{2}},\bze;0)\right)\cdot\bU^{(1)}_{\ep^2}(0,\bze)\right\}d{s_1}.
\ee
Here $\vV(\cdot,\cdot;0)$ denotes the orthogonal projection of
$\vV(\cdot,\cdot)$ onto $L^2_{-\infty,0}$. A similar notation has been used for $\nabla\chi_{\ep^2}$.

By incompressibility of  the velocity
\[
\bU^{(1)}_{\la} (t,\bx;0)=\vV(t,\bx;0)\cdot \chi_\la^{(1)}(t,\bx;0).
\]
Using  Cauchy-Schwarz inequality we
estimate
(\ref{1.5}) by
\be
\label{0303}
C_1(\bE|\bU^{(1)}_{\ep^2}(0,\bze)|^4)^{1/4}(\bE|\vV(0,\bze)|^4)^{1/4}\left(\int\limits_0^tds\int\limits_0^s\bE|\nabla\chi^{(1)}_{\ep^2}(\frac{{s_1}}{\ep^{2}},\bze;0)|^2d{s_1}\right)^{1/2}.
\ee
Since the fourth moment of $\bU^{(1)}_{\ep^2}$ remains bounded as
$\ep\da$, (\ref{0303})  can be estimated by using (\ref{1.10}) and
(\ref{40'}) as
\be
C_2\left[\int\limits_0^t\int\limits_0^K\frac{\ep^2[1-\exp(-\frac{2k^{2\bt}s}{\ep^2})]}{k^{2\bt}}\times{k^2\over
    k^{2\al-1}(\ep^2+k^{2\bt})^2}dkds\right]^{1/2},
\label{1.50}
\ee
which tends to $0$ for $\al+2\bt<2$,
in view of  the Dominant Convergence Theorem.

The second term in (\ref{31.12}) can be rewritten,
after a simple change of variables, as
\be
\label{223}
\frac{2}{\ep}\int\limits_0^tds\int\limits_0^sd{s_2}
\int\limits_{s_2}^{s}\bE\left\{\left(\vV(\frac{{s_2}}{\ep^{2}},
\bx_\ep({s_2}))\cdot\nabla\bU^{(1)}_{\ep^2}(\frac{{s_1}}{\ep^{2}},
\bx_\ep({s_2});\frac{s_2}{\ep^2})\right)
\cdot\bU^{(1)}_{\ep^2}(0,\bze)\right\}d{s_1}.
\ee
By the spectral representation of $\vV$  we obtain
\[
\nabla\bU^{(1)}_{\ep^2}(t,\bx;s)=
\int\limits_{R^d}\int\limits_{R^d}
\frac{\exp\left\{-{(|\bk|^{2\bt}+|\bk_1|^{2\bt})(t-s)\over
\ep^2}\right\}}{|\bk|^{2\bt}+\ep^2}\bW(s,\bx,d\bk,d\bk_1)
\]
where
\[
\bW(s,\bx,d\bk,d\bk_1):=\bk\cdot
\hat{\vV}_1(s,\bx,d\bk_1)\bk_1\otimes\hat{\vV}(s,\bx,d\bk)+
\bk\cdot
\hat{\vV}(s,\bx,d\bk_1)\bk\otimes\hat{\vV}_1(s,\bx,d\bk).
\]
Moreover by stationarity of the Lagrangian velocity experienced
by the particle (see \cite{port-stone}) and incompressibility of velocity
we have
\be
\label{0330}
{1\over \ep^2}\bE\left|\int\limits_{s_2}^{s}
\nabla\bU^{(1)}_{\ep^2}(\frac{s_1}{\ep^{2}},\bx_\ep({s_2});
\frac{s_2}{\ep^{2}})ds_1\right|^2=
{1\over \ep^2}\bE\left|\int\limits_{0}^{s-s_2}
\nabla\bU^{(1)}_{\ep^2}(\frac{{s_1}}{\ep^{2}},\bze;0)ds_1\right|^2=
\ee
\[
{1\over \ep^2}\bE\left|\int\limits_{R^d}\int\limits_{R^d}\frac{\ep^2\left[1-\exp\left\{-{(|\bk|^{2\bt}+|\bk_1|^{2\bt})(s-s_2)\over
    \ep^2}\right\}\right]}{(|\bk|^{2\bt}+\ep^2)
(|\bk|^{2\bt}+|\bk_1|^{2\bt})}\bW(0,\bze,d\bk,d\bk_1)\right|^2,
\]
which  equals
\[
\int\limits_{R^d}\int\limits_{R^d}\int\limits_{R^d}\int\limits_{R^d}\frac{\left[1-\exp\left\{-{(|\bk|^{2\bt}+|\bk_1|^{2\bt})(s-s_2)\over
    \ep^2}\right\}\right]}{(|\bk|^{2\bt}+\ep^2)(|\bk|^{2\bt}+|\bk_1|^{2\bt})}\times
\frac{\ep^2\left[1-\exp\left\{-{(|\bk'|^{2\bt}+|\bk'_1|^{2\bt})(s-s_2)\over
    \ep^2}\right\}\right]}{(|\bk'|^{2\bt}+|\bk'_1|^{2\bt})(|\bk'|^{2\bt}+\ep^2)}
\]
\[
\bE\left\{\bW(0,\bze,d\bk,d\bk_1)\bW(0,\bze,d\bk',d\bk_1')\right\}.
\]
By incompressibility of velocity again,
the above expectation equals
\[
\left[\left(|\bk_1|^2\mbox{tr}\hat{\bR}(\bk)\bk\cdot\hat{\bR}(\bk_1)
\bk+|\bk|^2\bk\cdot\hat{\bR}(\bk_1)\hat{\bR}(\bk)\bk_1\right)
\delta(\bk-\bk_1')\delta(\bk_1-\bk')\right.+
\]
\[
\left.\left(|\bk|^2\mbox{tr}\hat{\bR}(\bk)\bk\cdot\hat{\bR}(\bk_1)
\bk+|\bk_1|^2\bk\cdot\hat{\bR}(\bk_1)\hat{\bR}(\bk)\bk_1\right)
\delta(\bk-\bk')\delta(\bk_1-\bk_1')\right]
d\bk d\bk_1d\bk' d\bk_1'.
\]
Thus we can bound the utmost left hand side of (\ref{0330}) by
\be
\label{034}
C\left\{\int\limits_0^K\int\limits_0^K
\frac{k^3(k+k_1)}{(k^{2\bt}+\ep^2)^2k^{2\al-1}k_1^{2(\bt+\al)-1}}
\times {\ep^2\left[1-\exp\left\{-{ (k^{2\bt}+k_1^{2\bt})(s-s_2)
\over\ep^2}\right\}\right]\over(k^{2\bt}+k_1^{2\bt})}dk dk_1\right.+
\ee
\[
\left.\int\limits_0^K\int\limits_0^K{k^3(k+k_1)\over k^{2\al+\bt-1}k_1^{2\al+\bt-1}
(k^{2\bt}+\ep^2) (k^{2\bt}_1+\ep^2) }\times {\ep^2\left[1-\exp\left\{-{
(k^{2\bt}+k_1^{2\bt})(s-s_2)\over\ep^2}\right\}\right]\over(k^{2\bt}+k_1^{2\bt})}dk
dk_1\right\},
\]
with the constant $C$ independent of $\ep$. For $\al+2\bt<2$,
 (\ref{034})
tends to zero, as $\ep\da$, by the Dominant Convergence Theorem.
Thus, the second term of (\ref{31.12}) vanishes with $\ep$\kropa

\end{document}